\documentclass[11pt]{article}
\usepackage{bbm}
\usepackage{epsfig}
\usepackage{amssymb,amsmath,amscd, mdwmath}
\usepackage{graphicx,cite,url}

\begin{document}
\title{\bf Simulation of Wave Equation  on  Manifold using  DEC}

\author{
Zheng Xie$^1$\thanks{E-mail: lenozhengxie@yahoo.com.cn Tel./fax: +86
0739 5316081}~~~~Yujie Ma$^2$\thanks{ E-mail: yjma@mmrc.iss.ac.cn
This work is partially supported by CPSFFP (No. 20090460102), NKBRPC
(No. 2004CB318000) and NNSFC (No. 10871170) }
\\{\small
$1.$ Center of Mathematical Sciences, Zhejiang University
(310027),China}
\\ {\small $2.$ Key Laboratory of Mathematics Mechanization,}
\\ {\small  Chinese Academy of Sciences,  (100090), China}}

\date{}

\maketitle

\begin{abstract}

The classical numerical methods    play   important roles in solving
wave equation,  e.g. finite difference time domain method. However,
their computational domain are limited to flat space and the time.
This paper deals with the description of discrete exterior calculus
method for numerical simulation of wave equation. The advantage of
this method is that it can be used to compute equation on  the space
manifold and the time. The analysis of its stable condition and
error is also accomplished.

\end{abstract}


\vskip 0.2cm \noindent {\bf PASC(2010):}  43.20.+g, 02.30.Jr,
02.30.Mv, 02.40.Ky.

\section{Introduction}

The wave equation is the prototypical example of a hyperbolic
partial differential equation of waves, such as sound waves, light
waves and water waves. It arises in fields such as acoustics,
electromagnetism, and fluid dynamics\cite{morton,larsson}. To
investigate the predictions of wave equation of such phenomena it is
often necessary to approximate its solution numerically. A technique
suitable for providing numerical solutions to the  wave propagation
problem is the finite difference time domain (FDTD) method. This is
normally defined by looking for an approximate solution on a uniform
mesh of points and by replacing the derivatives in the differential
equation by difference quotients at   points of mesh. The
computational domain of this algorithm is limited to flat space and
the time \cite{bondeson,yee,bossavit,bossavit1,stern}.

Discrete exterior calculus (DEC) constitutes a discrete realization
of the exterior differential forms, and therefore, the right
framework in which to develop a  discretization for  differential
equations not just on flat space but   on  manifold \cite{
whitney,arnold,desbrun,leok,hiptmair,hyman,meyer,luo,xie-ye-ma}. The
differential operators   such as gradient, divergence, and Laplace
operator
 on manifold can also be naturally discretized  using  DEC.
 The numerical solution of
 wave equation on space manifold and the time by the methods of  DEC
 is obtained in this paper. For
this equation, an explicit scheme is derived. The analysis of this
scheme's stability  shows that the numerical solution becomes
unstable unless the time step is restricted.
\section{ DEC method for wave equation}

\subsection*{Wave equation}
 The wave equation is the prototypical
example of a hyperbolic partial differential equation. In its
simplest form, the wave equation refers to a scalar function $u$
that satisfies:
$$\frac{\partial^2 u}{\partial t^2}=c^2\Delta u,\eqno{(1)}$$
where $\Delta$ is the Laplace operator and  $c$ is the propagation
speed of the wave. More realistic differential equations for waves
allow   the speed of wave propagation to vary with the frequency of
the wave, a phenomenon known as dispersion. In this case, $c$ is
replaced by the phase velocity:
$$\frac{\partial^2 u}{\partial t^2}=\left(\frac{\omega}{k}\right)^2\Delta u.$$
Another common correction in realistic systems is that the speed is
depend on the amplitude of the wave, leading to a nonlinear wave
equation:
$$\frac{\partial^2 u}{\partial t^2}=c(u)^2\Delta u.$$

\subsection*{DEC scheme for wave equation}

  A discrete differential $k$-form, $k \in \mathbb{Z}$, is the
evaluation   of the differential $k$-form on all $k$-simplices. A
dual form is evaluated on the dual cell.  Suppose each simplex
contains its circumcenter. The circumcentric dual cell $D(\sigma_0)$
of simplex $\sigma_0$ is
 $$ D(\sigma_0):=\bigcup_{\sigma_0\in \sigma_1\in\cdots \in\sigma_r}
 \mathrm{Int}(c(\sigma_0)c(\sigma_1)\cdots c(\sigma_r) ),$$
where   $\sigma_i$ is all the simplices which contains
$\sigma_0$,..., $\sigma_{i-1}$, $c(\sigma_i)$ is the circumcenter of
$\sigma_i$. In DEC, the basic operators in differential geometry are
approximated as follows:
\begin{itemize}
  \item[1.] Discrete exterior derivative $d$, this operator is the transpose of the
incidence matrix of $k$-cells on $k+1$-cells.
  \item[2.] Discrete Hodge Star $\ast$, the operator     scales the cells by
  the volumes of the corresponding dual and
primal cells.
\item[3.] Discrete Laplace operator is   $  \ast^{-1}d^{T}\ast + d^{T}\ast
d.$
\end{itemize}

 For some situations, a source having azimuthal symmetry about
its axis is considered. In this case, the 2D triangular discrete
manifold as the space is only need to be considered.
 Now,  we show how to derive an explicit DEC  scheme for Eq.(1) in 2D   space manifold and the time.
 The wave equation in 3D
space and the time  can  also be computed by a similar approach.

 Take Fig.1  as an example for a part of 2D mesh, in
which $0$,..., $F$ are vertices,    $1$,...,$6$ are the
circumcenters of triangles,    $a$,...,$f$ are the circumcenters of
edges. Denote $l_{ij}$ as the length of line segment $(i,j)$ and
$A_{ijkl}$ as the area of  quadrangle $(i,j,k,l)$.
 $$
\begin{minipage}{0.99\textwidth}
\begin{center} \includegraphics[scale=0.3]{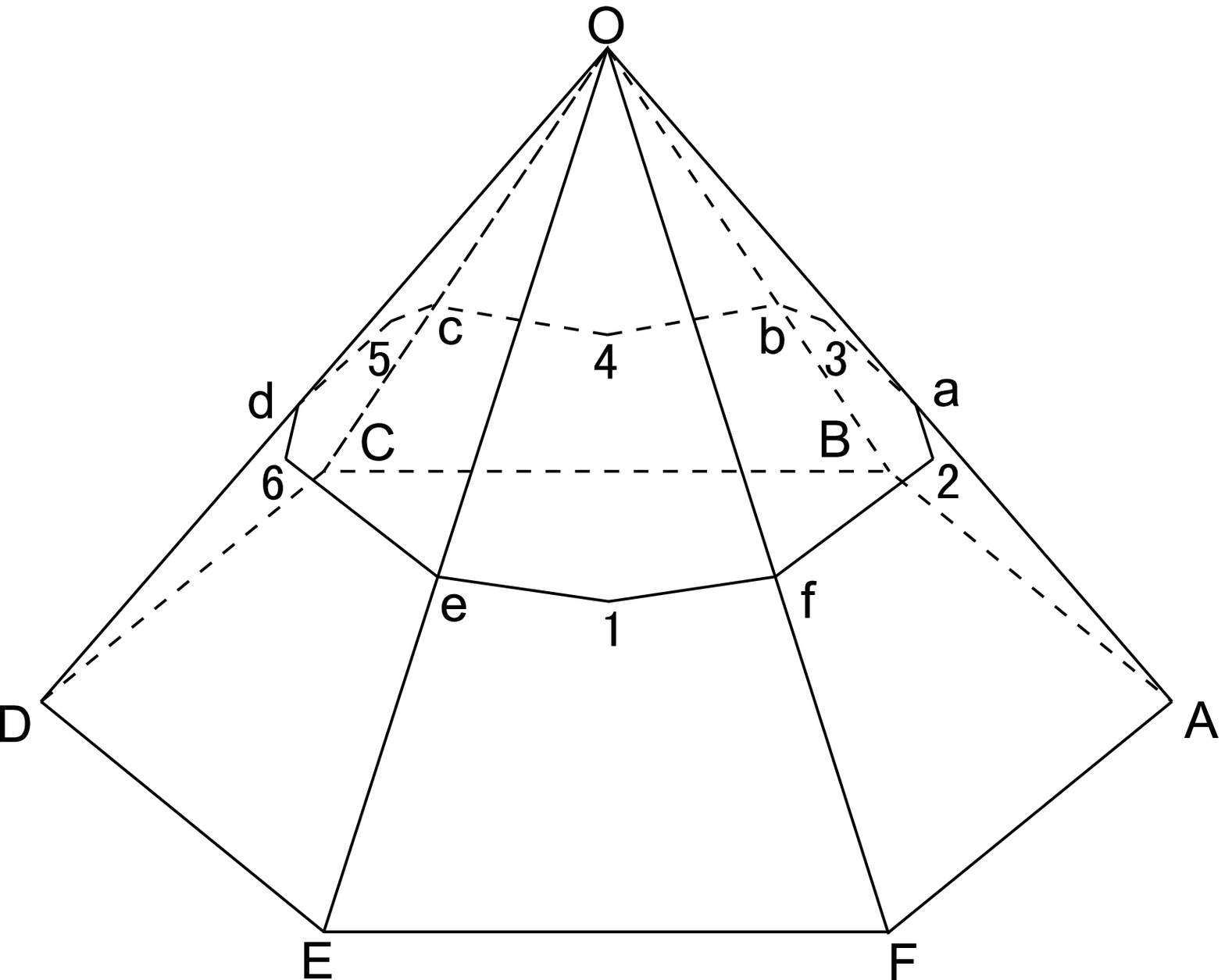}
\end{center}
\centering{ Fig.1. A part of 2D mesh}
\end{minipage}
$$Define $$l_{12}:=l_{1f}+l_{2f},~l_{23}:=l_{2a}+l_{3a},...,
   l_{61}:=l_{6e}+l_{1e},$$  and $$ P_{123456}:=
A_{01fe}+A_{02fa}+\cdots+A_{06de} .$$ The  diffusion term $\Delta u$
at vertice $0$ is approximated using discrete Laplace operator as
follows:
$$ \begin{array}{lll}\Delta u_0&\approx& \dfrac{1}{P_{123456}} \left(\dfrac{l_{23}}{l_{A0}}
(u_A-u_0)+\dfrac{l_{34}}{l_{B0}}(u_B-u_0)+\dfrac{l_{45}}{l_{C0}}(u_C-u_0)\right.\\
&&\left.+\dfrac{l_{56}}{l_{D0}}(u_D-u_0)+\dfrac{l_{16}}{l_{E0}}(u_E-u_0)+\dfrac{l_{12}}{l_{F0}}(u_F-u_0)
\right).\end{array}\eqno{(2)}
$$
The temporal derivative  is approximated by middle time differences
as follows:
$$\frac{\partial^2u^n}{\partial t^2}\approx \frac{1}{(\Delta t)^2}(u^{n+1}-2u^n+u^{n-1}),\eqno{(3)}$$
where $\Delta t$ is  uniform spacing, and $n\Delta t$ is the
coordinate of time. The approximation of Eq.(1) generated by
substituting the left-hand sides of (2) and (3) into (1), thus
satisfies
$$\begin{array}{lll} \mathrm{{Right}}(2)^{n-1} &=&\dfrac{1}{(c\Delta
t)^2}\left(u^{n}_0-2u^{n-1}_0+u^{n-2}_0\right).
\end{array}\eqno{(4)}$$

\section{Stability, convergence  and accuracy}

\subsection*{Stability}
The Courant-Friedrichs-Lewy condition   is a necessary condition for
stability while solving certain partial differential equations
numerically.    Now, this condition is derived for scheme (4).
First, this  DEC scheme is  decomposed into temporal and spacial
eigenvalue problems.

The temporal eigenvalue problem:
$$\dfrac{\partial^2u^n}{\partial t^2}=\Lambda u^n$$
It can be approximated by difference equation
$$\dfrac{u^{n+1}_0-2u^{n}_0+u^{n-1}_0}{(\Delta t)^2}=\Lambda u^n_0.\eqno{(5)}$$
Supposing $$u^{n+1}_0= u^n_0 \cos( \Delta t)~~~~u^{n-1}_0=u^{n }_0
\cos( -\Delta t)$$  and substituting those into Eq.(5), we obtain
$$\dfrac{\cos( \Delta t)+\cos( -\Delta t)-2}{ (\Delta t)^2}  =\Lambda,$$
therefore
$$-\dfrac{ 4}{ ( \Delta t)^2}  \leq \Lambda\leq 0.$$
This is the   stable condition for the temporal eigenvalue problem.

The spacial eigenvalue problem:
$$c^2\Delta u=\Lambda u$$
It can be approximated by difference equation (6) based on Fig.1.
  $$
\begin{array}{lll}
   \dfrac{P_{123456 }}{c^2}\Lambda u_0 &=& \dfrac{l_{23} }{l_{A0}}(u_A
-u_0 )+\dfrac{l_{34} }{l_{B0} }(u_B
-u_0 )+\dfrac{l_{45} }{l_{C0} }(u_C -u_0 )\\
&&+\dfrac{l_{56} }{l_{D0} }(u_D -u_0 )+\dfrac{l_{16} }{l_{E0} }(u_E
-u_0 ) +\dfrac{l_{12} }{l_{F0} }(u_F -u_0)
\end{array}\eqno{(6)}
$$
Let $ u_i= u_0 \cos(c l_{0i}) $ and substitute into Eq.(6) to obtain
  $$\begin{array}{lll}
  \dfrac{P_{123456}}{c^2}\Lambda &=& \dfrac{l_{23} }{l_{A0}}(\cos(c l_{0A})
-1 )+\dfrac{l_{34} }{l_{B0} }(\cos(c l_{0B}) -1 )+\dfrac{l_{45}
}{l_{C0} }(\cos(c l_{0C}) -1 )\\
&&+\dfrac{l_{56} }{l_{D0} }(\cos(c l_{0D}) -1 )+\dfrac{l_{16}
}{l_{E0} }(\cos(c l_{0E}) -1 ) +\dfrac{l_{12} }{l_{F0} }(\cos(c
l_{0F}) -1 )
\end{array}$$So we have
$$-\dfrac{2c^2}{P_{123456}}\left(\dfrac{l_{23} }{l_{A0} } +\dfrac{l_{34}
}{l_{B0} } +\dfrac{l_{45} }{l_{C0} } +\dfrac{l_{56} }{l_{D0} }
+\dfrac{l_{16} }{l_{E0} } +\dfrac{l_{12} }{l_{F0}  }\right)
\leq\Lambda\leq 0.$$ In order to keep the stability of scheme (4),
we need
$$-\dfrac{ 4}{ ( \Delta t)^2}\leq -\dfrac{2c^2}{P_{123456}}\left(\dfrac{l_{23} }{l_{A0} } +\dfrac{l_{34}
}{l_{B0} } +\dfrac{l_{45} }{l_{C0} } +\dfrac{l_{56} }{l_{D0} }
+\dfrac{l_{16} }{l_{E0} } +\dfrac{l_{12} }{l_{F0} }\right)
\eqno{(7)}$$for all vertices, namely
$$c \Delta t \leq {\mathrm{Min}}_{0\in V} \sqrt{\dfrac{2P_{123456}}{{\dfrac{l_{23} }{l_{A0} } +\dfrac{l_{34}
}{l_{B0} } +\dfrac{l_{45} }{l_{C0} } +\dfrac{l_{56} }{l_{D0} }
+\dfrac{l_{16} }{l_{E0} } +\dfrac{l_{12} }{l_{F0} }}}}.$$

\subsection*{Convergence}

By the definition of truncation error, the solution $\tilde{u}$ of
the Eq.(1) satisfies the same relation as scheme (4) except for an
additional term $O(\Delta t)^2$ on the right hand side. Thus the
error $X^n_i=\tilde{u}^n_i-u^n_i$ is determined from the relation
$$\begin{array}{lll} X^n_0 &=&
2X^{n-1}_0- X^{n-2}_0+
  \dfrac{(c\Delta t)^2}{P_{123456}}\left(\dfrac{l_{23}}{l_{A0} }(X^{n-1}_A-X^{n-1}_0)\right.\\
& &+\left. \dfrac{l_{34} }{l_{B0} }(X^{n-1}_B -X^{n-1}_0
)+\dfrac{l_{45} }{l_{C0} }(X^{n-1}_C-X^{n-1}_0) +\dfrac{l_{56}
}{l_{D0} }(X^{n-1}_D-X^{n-1}_0)\right.
\\
& &+\left.\dfrac{l_{16} }{l_{E0} }(X^{n-1}_E-X^{n-1}_0)
+\dfrac{l_{12} }{l_{F0} }(X^{n-1}_F-X^{n-1}_0) \right)
             + O(\Delta t)^2.
\end{array} \eqno{(8)}$$
 Define   $$|X^n|=\mathrm{Max}_{i\in V} |X^n_i|.$$ From condition (7), we have
 $$ \dfrac{(c\Delta t)^2}{P_{123456}}\left({\dfrac{l_{23} }{l_{A0} } +\dfrac{l_{34} }{l_{B0} }
+\dfrac{l_{45} }{l_{C0} } +\dfrac{l_{56} }{l_{D0} } +\dfrac{l_{16}
}{l_{E0} } +\dfrac{l_{12} }{l_{F0} }} \right) <2.$$ Hence the
coefficient of $X^{n-1}_0$ in Eq.(8) is nonnegative. It follows that
$$\begin{array}{lll}|X^{n}_0|&\leq  &  2|X^{n-1}|+ |X^{n-2}|+O(\Delta t)^2,
\end{array} $$ and hence that
$$|X^{n}| \leq   2|X^{n-1}|+ |X^{n-2}|+O(\Delta t)^2.$$ Iterating $n$, we
obtain
$$ |X^n|<M_1  |X^1|+M_0 |X^0|+ O(\Delta t)^2,$$ where $M_1$, $M_0$ are
finite value defined on $n$. Since the initial conditions ensure
$X^0=0$ and $X^1=0$, we have
$$\lim_{\Delta t\rightarrow 0 } |X^n|=0.$$
That is to say the numerical solution approaches the exact solution
as the step size goes to $0$, and scheme (4) is   convergent.

\subsection*{Accuracy}

In scheme (4), the space derivative of  is approximated by first
order difference. Equivalently, $u$ is approximated by linear
interpolation functions. Consulting the
  definition about accuracy of finite volume method, we
can   say that scheme (4) has first order spacial accuracy. Scheme
(4)  has second order temporal accuracy, and second order spacial
accuracy on rectangular grid with uniform spacing.

\section{Algorithm Implementation}

 The
implementation of DEC scheme (4) for wave  equation consists of the
following steps:
\begin{itemize}
  \item [1.]Set the simulation parameters. These
are the dimensions of the computational mesh and the size of the
time step, etc.;

 \item [2.] Initialize the mesh indexes.

  \item [3.]Assign current transmitted signal.

  \item [4.]Compute the value of all spatial nodes  and temporarily store the
  result  in the circular buffer for further computation.

  \item [5.]Visualize the currently computed grid of spatial nodes.

   \item [6.]Repeat the  process  from the step 3, until reach the desired total number
   of iterations.
\end{itemize}
The flowchart of the   scheme (4) can be seen in Fig.2.
$$
\begin{minipage}{0.99\textwidth}
\begin{center} \includegraphics[scale=0.35]{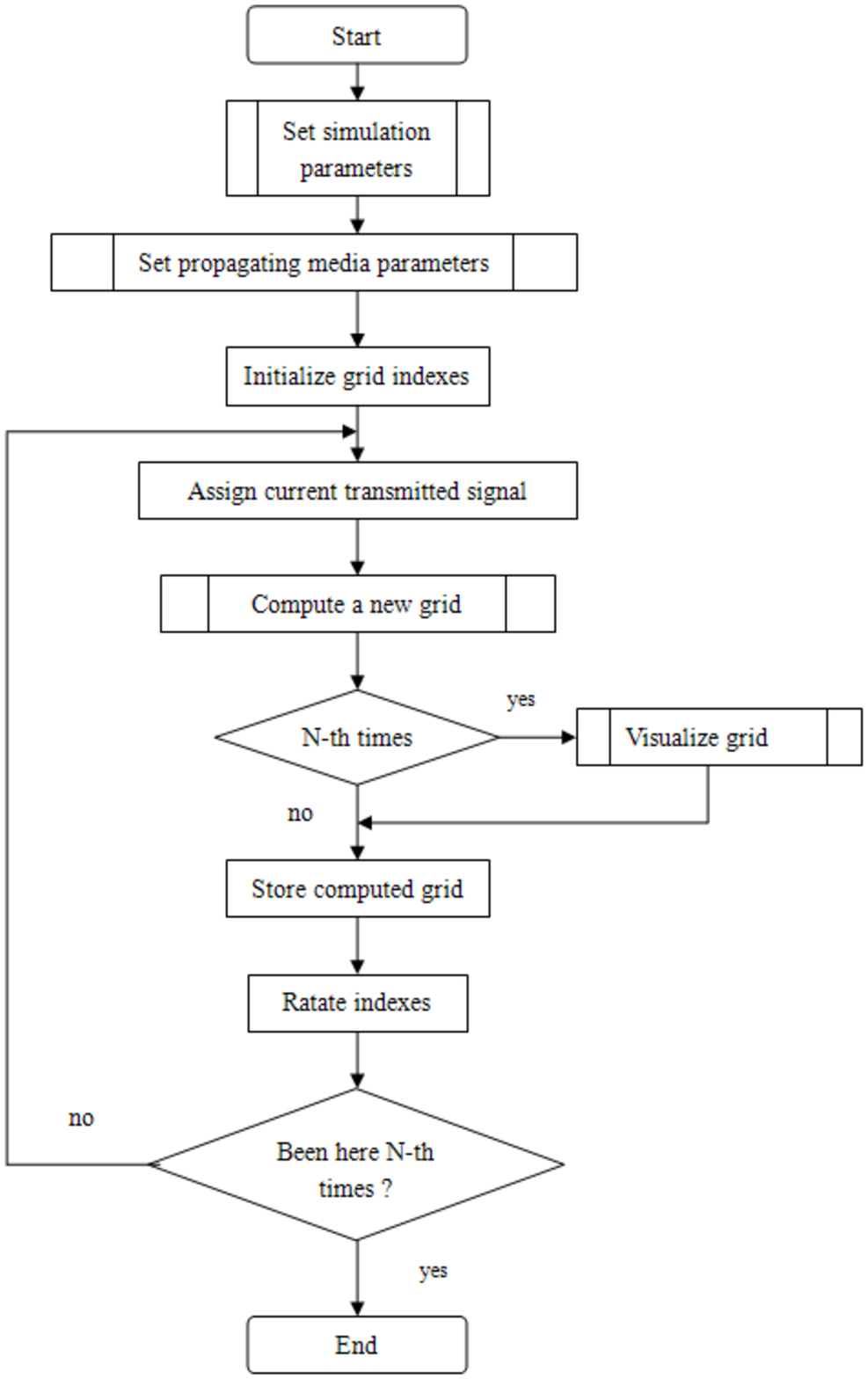}
\end{center}
\centering{ Fig.2. The flowchart of  scheme (4)}
\end{minipage}
$$

The Fig.3 and Fig.4 show the numerical simulation of Gaussian pulse
propagating on  the sphere and rabbit by scheme (4)  on C\#
platform.
$$
\begin{minipage}{0.99\textwidth}
\begin{center} \includegraphics[scale=0.62]{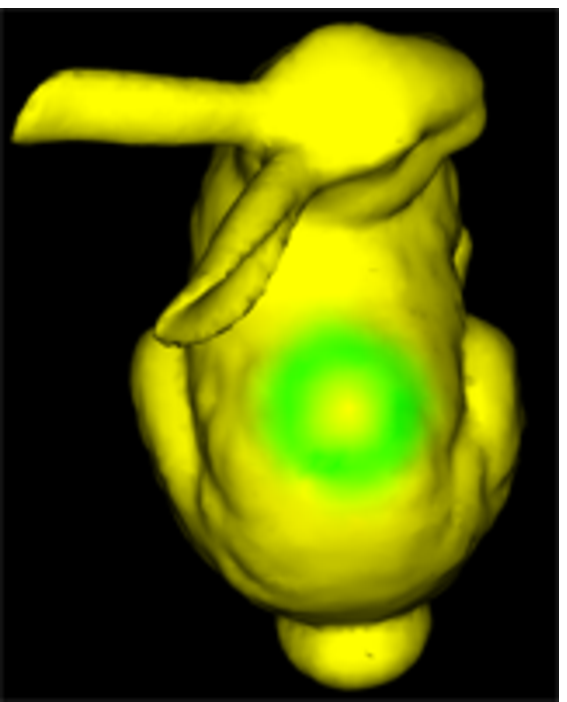}
 \includegraphics[scale=0.62]{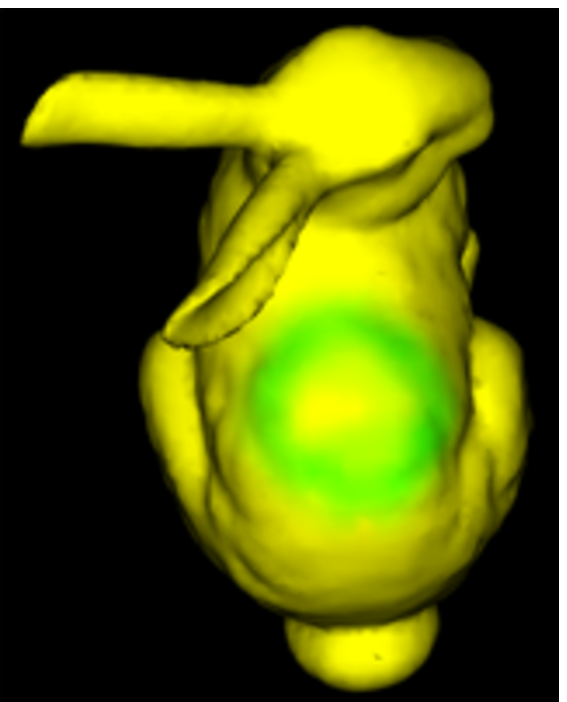}
 \includegraphics[scale=0.62]{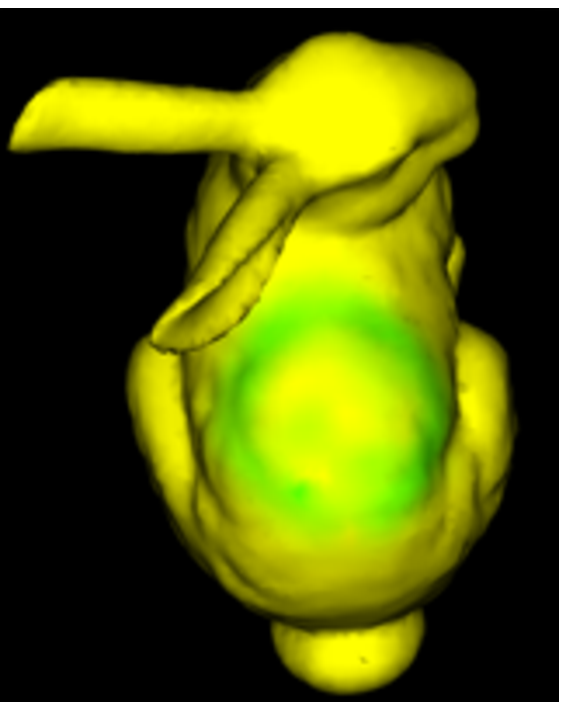}
\end{center}
\centering{ Fig.3. The propagation of Gaussian pulse on a rabbit}
\end{minipage}
$$
$$
\begin{minipage}{0.99\textwidth}
\begin{center} \includegraphics[scale=0.5]{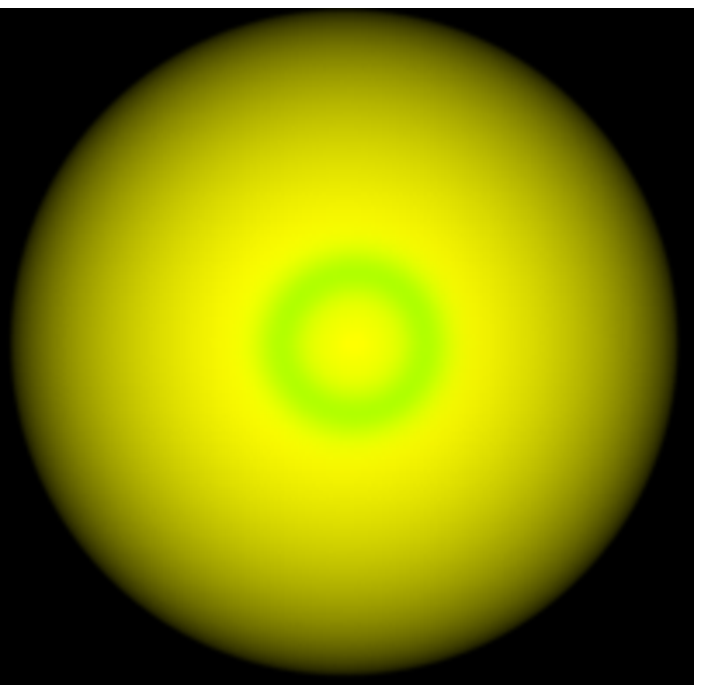}
 \includegraphics[scale=0.5]{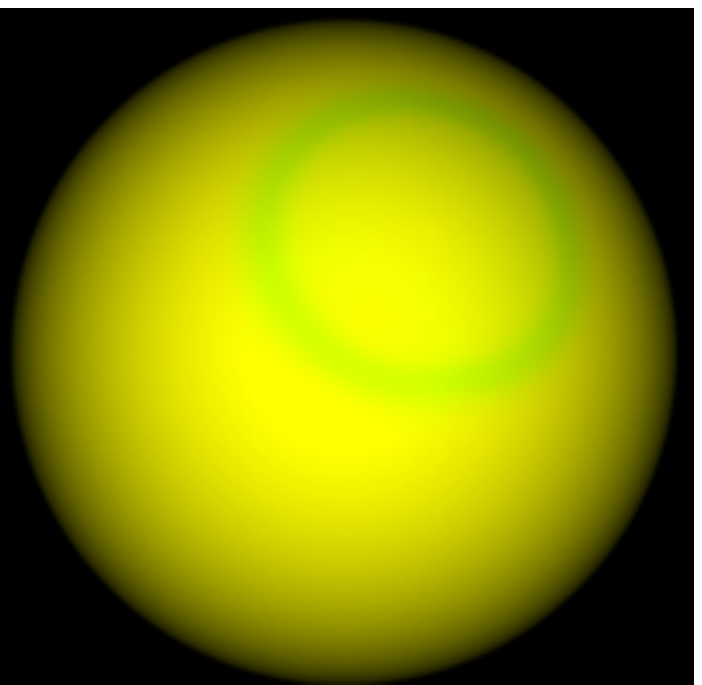}
 \includegraphics[scale=0.5]{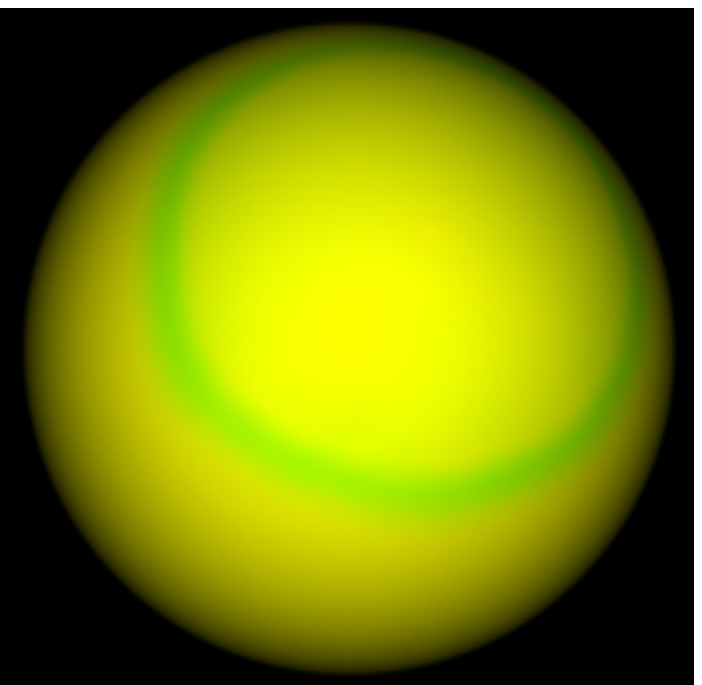}
\end{center}
\centering{ Fig.4. The propagation of Gaussian pulse on a sphere}
\end{minipage}
$$

\section{Discussion}
The DEC scheme for Laplace operator here can also be used to
simulate the heat equation, Laplace equation and Poisson equation on
manifold.

\subsection*{Discrete Laplace equation}

The discrete Laplace  equation  on surface of regular tetrahedron
(Fig.(5)) is
$$\left(
                   \begin{array}{cccc}
                     1 & 1 & 1 & -3 \\
                     1 & 1 & -3 & 1 \\
                     1 & -3 & 1 & 1 \\
                     -3 & 1 & 1 &  1 \\
                   \end{array}
                 \right)\left(
                          \begin{array}{c}
                            u_A \\
                            u_B \\
                            u_C \\
                            u_D \\
                          \end{array}
                        \right)=\left(
                                  \begin{array}{c}
                                    0 \\
                                    0 \\
                                    0 \\
                                    0 \\
                                  \end{array}
                                \right)\eqno{(9)}
$$
The solution of Eq.(9) is $$u_A= u_B=u_C=u_D=C ,$$ where $C$ is
arbitrary constant. Eq.(9) is an imprecise approximation of
Laplace's equations on a sphere. Obviously, this equation has
constant solution.
$$
\begin{minipage}{0.99\textwidth}
\begin{center} \includegraphics[scale=0.3]{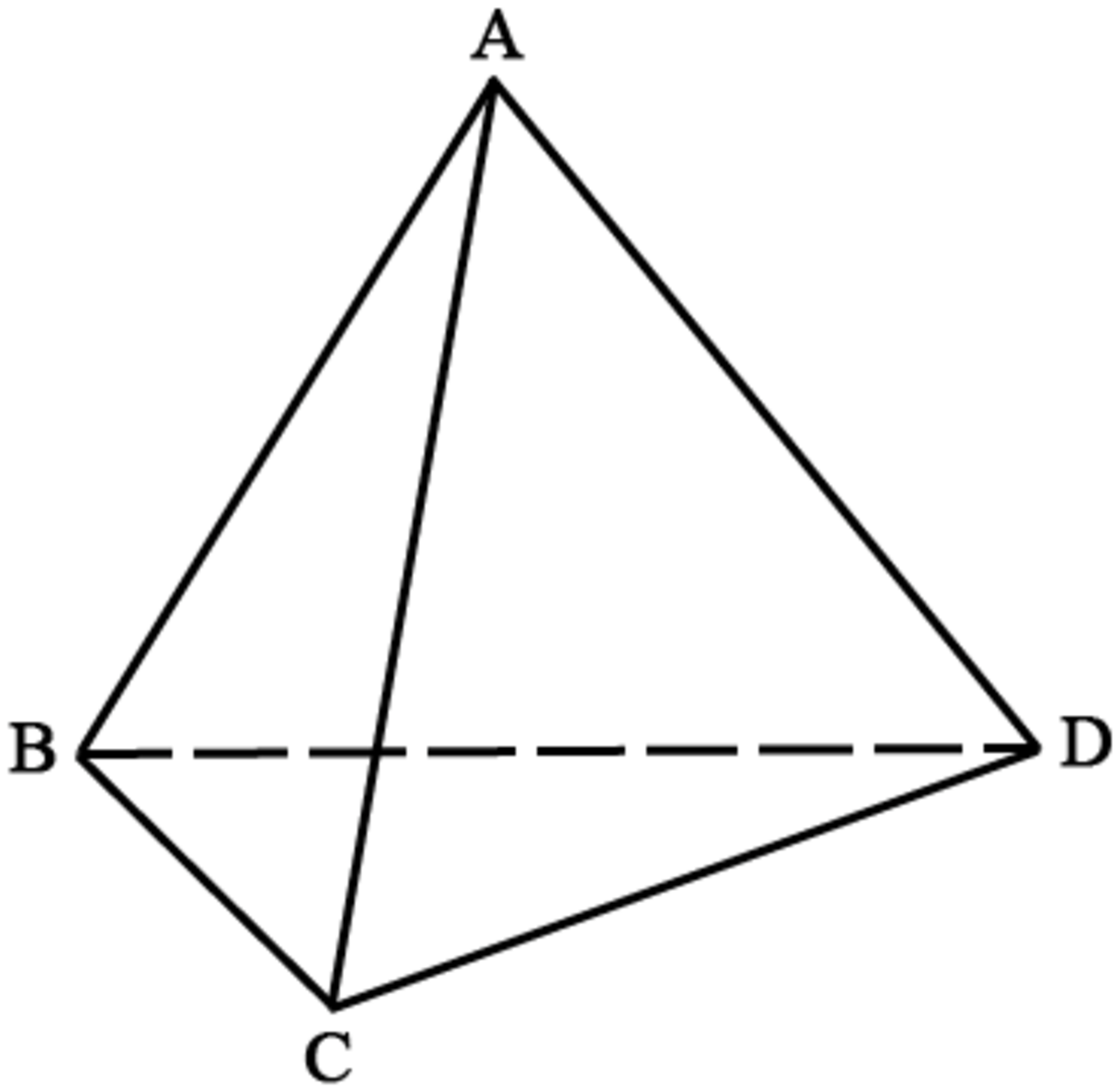}
\end{center}
\centering{ Fig.5. The surface of regular tetrahedron}
\end{minipage}
$$

\subsection*{Discrete Poisson equation}
 Consider a discrete Poisson equation on surface of regular tetrahedron.
Suppose the boundary condition is $$u_A=H ,$$ then discrete Poisson
equation on Fig.(5) is
$$\left(
                   \begin{array}{cccc}
                     3 & -1 & -1   \\
                     -1 & 3 & -1 \\
                     -1 & -1 & 3  \\
                   \end{array}
                 \right)\left(
                          \begin{array}{c}
                            u_B \\
                            u_C \\
                            u_D \\
                          \end{array}
                        \right)=\left(
                                  \begin{array}{c}
                                    H \\
                                    H \\
                                    H \\
                                  \end{array}
                                \right)\eqno{(10)}
$$
The solution of Eqs.(10) is $$u_B=u_C=u_D=H .$$

\subsection*{Discrete heat equation}
The heat equation of temperature $u$ is
$$\frac{\partial u}{\partial t}=c\Delta u,$$ which can be
approximated as
$$ u^{n}_0  =u^{n-1}_0+c\Delta
t~ \mathrm{{Right}}(2)^{n-1}.\eqno{(11)}
 $$
  The Fig.6 shows
the heat diffusion  of a constant heat source on the sphere
simulated by scheme (11).
$$
\begin{minipage}{0.99\textwidth}
\begin{center} \includegraphics[scale=0.5]{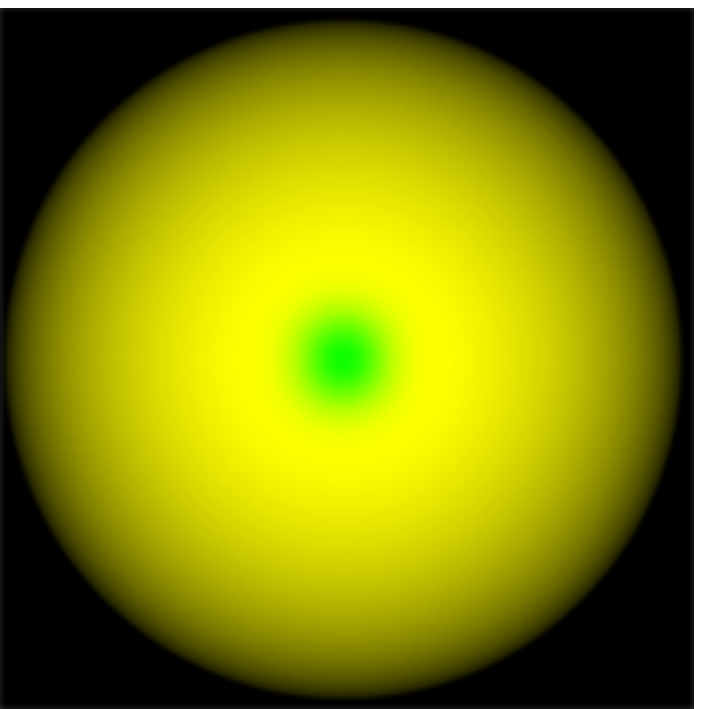}
 \includegraphics[scale=0.52]{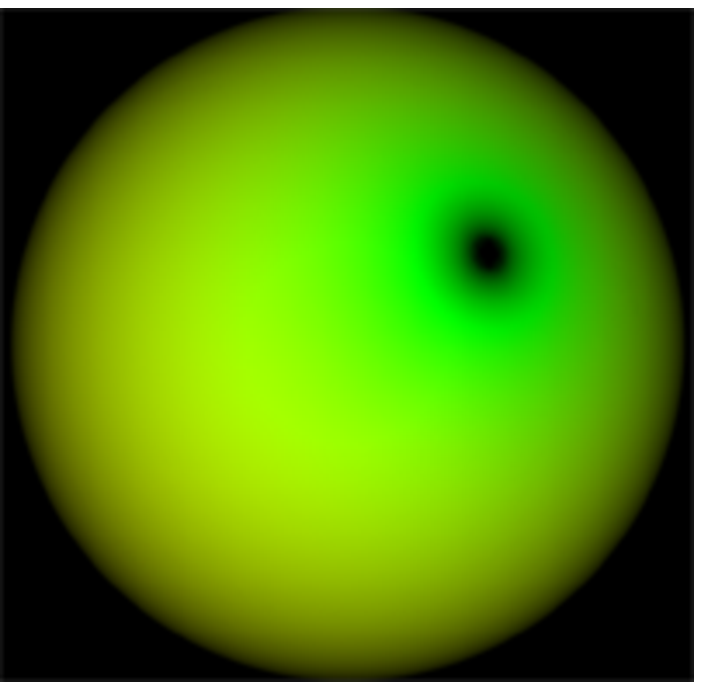}
 \includegraphics[scale=0.5]{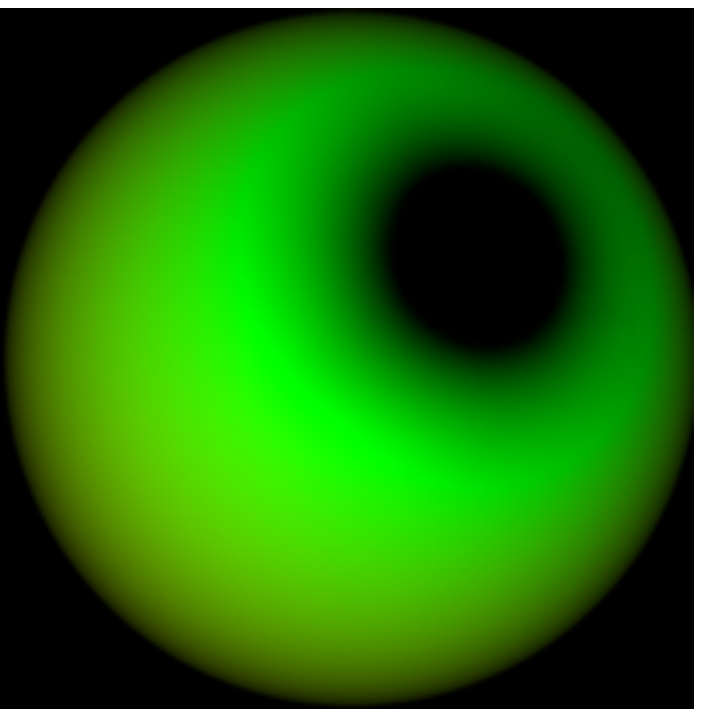}
\end{center}
\centering{ Fig.6. The heat diffusion on a sphere}
\end{minipage}
$$

\end{document}